# ON THE ADAPTIVE ELASTIC-NET WITH A DIVERGING NUMBER OF PARAMETERS

By Hui Zou[1] and Hao Helen Zhang[2]

*University of Minnesota and North Carolina State University*

We consider the problem of model selection and estimation in situations where the number of parameters diverges with the sample size. When the dimension is high, an ideal method should have the oracle property [*J. Amer. Statist. Assoc.* **96** (2001) 1348–1360] and [*Ann. Statist.* **32** (2004) 928–961] which ensures the optimal large sample performance. Furthermore, the high-dimensionality often induces the collinearity problem, which should be properly handled by the ideal method. Many existing variable selection methods fail to achieve both goals simultaneously. In this paper, we propose the adaptive elastic-net that combines the strengths of the quadratic regularization and the adaptively weighted lasso shrinkage. Under weak regularity conditions, we establish the oracle property of the adaptive elastic-net. We show by simulations that the adaptive elastic-net deals with the collinearity problem better than the other oracle-like methods, thus enjoying much improved finite sample performance.

**1. Introduction.**

1.1. *Background.* Consider the problem of model selection and estimation in the classical linear regression model

$$\mathbf{y} = \mathbf{X}\boldsymbol{\beta}^* + \boldsymbol{\varepsilon}, \tag{1.1}$$

where $\mathbf{y} = (y_1, \ldots, y_n)^T$ is the response vector and $\mathbf{x}_j = (x_{1j}, \ldots, x_{nj})^T, j = 1, \ldots, p$, are the linearly independent predictors. Let $\mathbf{X} = [\mathbf{x}_1, \ldots, \mathbf{x}_p]$ be the

Received December 2007; revised March 2008.
[1]Supported by NSF Grant DMS-07-06733.
[2]Supported by NSF Grant DMS-06-45293 and National Institutes of Health Grant NIH/NCI R01 CA-085848.

predictor matrix. Without loss of generality, we assume the data are centered, so the intercept is not included in the regression function. Throughout this paper, we assume the errors are identically and independently distributed with zero mean and finite variance $\sigma^2$. We are interested in the sparse modeling problem where the true model has a sparse representation (i.e., some components of $\boldsymbol{\beta}^*$ are exactly zero). Let $\mathcal{A} = \{j : \beta_j^* \neq 0, j = 1, 2, \ldots, p\}$. In this work, we call the size of $\mathcal{A}$ the intrinsic dimension of the underlying model. We wish to discover the set $\mathcal{A}$ and estimate the corresponding coefficients.

Variable selection is fundamentally important for knowledge discovery with high-dimensional data [Fan and Li (2006)] and it could greatly enhance the prediction performance of the fitted model. Traditional model selection procedures follow best-subset selection and its step-wise variants. However, best-subset selection is computationally prohibitive when the number of predictors is large. Furthermore, as analyzed by Breiman (1996), subset selection is unstable; thus, the resulting model has poor prediction accuracy. To overcome the fundamental drawbacks of subset selection, statisticians have recently proposed various penalization methods to perform simultaneous model selection and estimation. In particular, the lasso [Tibshirani (1996)] and the SCAD [Fan and Li (2001)] are two very popular methods due to their good computational and statistical properties. Efron et al. (2004) proposed the LARS algorithm for computing the entire lasso solution path. Knight and Fu (2000) studied the asymptotic properties of the lasso. Fan and Li (2001) showed that the SCAD enjoys the oracle property, that is, the SCAD estimator can perform as well as the oracle if the penalization parameter is appropriately chosen.

1.2. *Two fundamental issues with the $\ell_1$ penalty.* The lasso estimator [Tibshirani (1996)] is obtained by solving the $\ell_1$ penalized least squares problem

$$\widehat{\boldsymbol{\beta}}(\text{lasso}) = \arg\min_{\boldsymbol{\beta}} \|\mathbf{y} - \mathbf{X}\boldsymbol{\beta}\|_2^2 + \lambda\|\boldsymbol{\beta}\|_1, \tag{1.2}$$

where $\|\boldsymbol{\beta}\|_1 = \sum_{j=1}^p |\beta_j|$ is the $\ell_1$-norm of $\boldsymbol{\beta}$. The $\ell_1$ penalty enables the lasso to simultaneously regularize the least squares fit and shrink some components of $\widehat{\boldsymbol{\beta}}(\text{lasso})$ to zero for some appropriately chosen $\lambda$. The entire lasso solution paths can be computed by the LARS algorithm [Efron et al. (2004)]. These nice properties make the lasso a very popular variable selection method.

Despite its popularity, the lasso does have two serious drawbacks: namely, the lack of oracle property and instability with high-dimensional data. First of all, the lasso does not have the oracle property. Fan and Li (2001) first pointed out that asymptotically the lasso has nonignorable bias for estimating the nonzero coefficients. They further conjectured that the lasso may not have the oracle property because of the bias problem. This conjecture was recently proven in Zou (2006). Zou (2006) further showed that the lasso could be inconsistent for model selection unless the predictor matrix (or the design matrix) satisfies a rather strong condition. Zou (2006) proposed the following adaptive lasso estimator

$$(1.3) \qquad \widehat{\boldsymbol{\beta}}(\text{AdaLasso}) = \arg\min_{\boldsymbol{\beta}} \|\mathbf{y} - \mathbf{X}\boldsymbol{\beta}\|_2^2 + \lambda \sum_{j=1}^p \hat{w}_j |\beta_j|,$$

where $\{\hat{w}_j\}_{j=1}^p$ are the adaptive data-driven weights and can be computed by $\hat{w}_j = (|\hat{\beta}_j^{ini}|)^{-\gamma}$, where $\gamma$ is a positive constant and $\widehat{\boldsymbol{\beta}}^{ini}$ is an initial root-$n$ consistent estimate of $\boldsymbol{\beta}$. Zou (2006) showed that, with an appropriately chosen $\lambda$, the adaptive lasso performs as well as the oracle. Candes, Wakin and Boyd (2008) used the adaptive lasso idea to enhance sparsity in sparse signal recovery via the reweighted $\ell_1$ minimization.

Secondly, the $\ell_1$ penalization methods can have very poor performance when there are highly correlated variables in the predictor set. The collinearity problem is often encountered in high-dimensional data analysis. Even when the predictors are independent, as long as the dimension is high, the maximum sample correlation can be large, as shown in Fan and Lv (2008). Collinearity can severely degrade the performance of the lasso. As shown in Zou and Hastie (2005), the lasso solution paths are unstable when predictors are highly correlated. Zou and Hastie (2005) proposed the elastic-net as an improved version of the lasso for analyzing high-dimensional data. The elastic-net estimator is defined as follows:

$$(1.4) \quad \widehat{\boldsymbol{\beta}}(\text{enet}) = \left(1 + \frac{\lambda_2}{n}\right)\left\{\arg\min_{\boldsymbol{\beta}} \|\mathbf{y} - \mathbf{X}\boldsymbol{\beta}\|_2^2 + \lambda_2 \|\boldsymbol{\beta}\|_2^2 + \lambda_1 \|\boldsymbol{\beta}\|_1\right\}.$$

If the predictors are standardized (each variable has mean zero and $L_2$-norm one), then we should change $(1 + \frac{\lambda_2}{n})$ to $(1 + \lambda_2)$ as in Zou and Hastie (2005). The $\ell_1$ part of the elastic-net performs automatic variable selection, while the $\ell_2$ part stabilizes the solution paths and, hence, improves the prediction. In an orthogonal design where the lasso is shown to be optimal Donoho et al. (1995), the elastic-net automatically reduces to the lasso. However, when the correlations among the predictors become high, the elastic-net can significantly improve the prediction accuracy of the lasso.

1.3. *The adaptive elastic-net.* The adaptively weighted $\ell_1$ penalty and the elastic-net penalty improve the lasso in two different directions. The adaptive lasso achieves the oracle property of the SCAD and the elastic-net handles the collinearity. However, following the arguments in Zou and Hastie (2005) and Zou (2006), we can easily see that the adaptive lasso inherits the instability of the lasso for high-dimensional data, while the elastic-net lacks the oracle property. Thus, it is natural to consider combining the ideas of the adaptively weighted $\ell_1$ penalty and the elastic-net regularization to obtain a better method that can improve the lasso in both directions. To this end, we propose the adaptive elastic-net that penalizes the squared error loss using a combination of the $\ell_2$ penalty and the adaptive $\ell_1$ penalty. Since the adaptive elastic-net is designed for high-dimensional data analysis, we study its asymptotic properties under the assumption that the dimension diverges with the sample size.

Pioneering papers on asymptotic theories with diverging number of parameters include [Huber (1988) and Portnoy (1984)] which studied the $M$-estimators. Recently, Fan, Peng and Huang (2005) studied a semi-parametric model with a growing number of nuisance parameters, whereas Lam and Fan (2008) investigated the profile likelihood ratio inference for the growing number of parameters. In particular, our work is influenced by Fan and Peng (2004) who studied the oracle property of nonconcave penalized likelihood estimators. Fan and Peng (2004) provocatively argued that it is important to study the validity of the oracle property when the dimension diverges. We would like to know whether the adaptive elastic-net enjoys the oracle property with a diverging number of predictors. This question will be thoroughly investigated in this paper.

The rest of the paper is organized as follows. In Section 2, we introduce the adaptive elastic-net. Statistical theory, including the oracle property, of the adaptive elastic-net is established in Section 3. In Section 4, we use simulation to compare the finite sample performance of the adaptive elastic-net with the SCAD and other competitors. Section 5 discusses how to combine SIS of Fan and Lv (2008) and the adaptive elastic-net to deal with the ultra-high dimension cases. Technical proofs are presented in Section 6.

**2. Method.** The adaptive elastic-net can be viewed as a combination of the elastic-net and the adaptive lasso. Suppose we first compute the elastic-net estimator $\widehat{\boldsymbol{\beta}}(\text{enet})$ as defined in (1.4), and then we construct the adaptive weights by

$$(2.1) \qquad \hat{w}_j = (|\hat{\beta}_j(\text{enet})|)^{-\gamma}, \qquad j = 1, 2, \ldots, p,$$

where $\gamma$ is a positive constant. Now we solve the following optimization problem to get the adaptive elastic-net estimates

$$(2.2) \quad \begin{aligned} &\widehat{\boldsymbol{\beta}}(\text{AdaEnet}) \\ &= \left(1 + \frac{\lambda_2}{n}\right)\left\{\arg\min_{\boldsymbol{\beta}} \|\mathbf{y} - \mathbf{X}\boldsymbol{\beta}\|_2^2 + \lambda_2\|\boldsymbol{\beta}\|_2^2 + \lambda_1^* \sum_{j=1}^p \hat{w}_j |\beta_j|\right\}. \end{aligned}$$

From now on, we write $\widehat{\boldsymbol{\beta}} = \widehat{\boldsymbol{\beta}}(\text{AdaEnet})$ for the sake of convenience.

If we force $\lambda_2$ to be zero in (2.2), then the adaptive elastic-net reduces to the adaptive lasso. Following the arguments in Zou and Hastie (2005), we can easily show that in an orthogonal design the adaptive elastic-net reduces to the adaptive lasso, regardless the value of $\lambda_2$. This is desirable because, in that setting, the adaptive lasso achieves the optimal minimax risk bound [Zou (2006)]. The role of the $\ell_2$ penalty in (2.2) is to further regularize the adaptive lasso fit whenever the collinearity may cause serious trouble.

We know the elastic-net naturally adopts a sparse representation. One can use $\hat{w}_j = (|\hat{\beta}_j(\text{enet})| + 1/n)^{-\gamma}$ to avoid dividing zeros. We can also define $\hat{w}_j = \infty$ when $\hat{\beta}_j(\text{enet}) = 0$. Let $\widehat{\mathcal{A}}_{\text{enet}} = \{j : \hat{\beta}_j(\text{enet}) \neq 0\}$ and $\widehat{\mathcal{A}}_{\text{enet}}^c$ denotes its complement set. Then, we have $\widehat{\boldsymbol{\beta}}_{\widehat{\mathcal{A}}_{\text{enet}}^c} = 0$ and

$$
\begin{aligned}
(2.3) \quad \widehat{\boldsymbol{\beta}}_{\widehat{\mathcal{A}}_{\text{enet}}} &= \left(1 + \frac{\lambda_2}{n}\right) \\
&\quad \times \left\{ \arg\min_{\boldsymbol{\beta}} \|\mathbf{y} - \mathbf{X}_{\widehat{\mathcal{A}}_{\text{enet}}} \boldsymbol{\beta}\|_2^2 + \lambda_2 \|\boldsymbol{\beta}\|_2^2 + \lambda_1^* \sum_{j \in \widehat{\mathcal{A}}_{\text{enet}}} \hat{w}_j |\beta_j| \right\},
\end{aligned}
$$

where $\boldsymbol{\beta}$ in (2.3) is a vector of length $|\widehat{\mathcal{A}}_{\text{enet}}|$, the size of $\widehat{\mathcal{A}}_{\text{enet}}$.

The $\ell_1$ regularization parameters $\lambda_1^*$ and $\lambda_1$ are directly responsible for the sparsity of the estimates. Their values are allowed to be different. On the other hand, we use the same $\lambda_2$ for the $\ell_2$ penalty component in the elastic-net and the adaptive elastic-net estimators, because the $\ell_2$ penalty offers the same kind of contribution in both estimators.

**3. Statistical theory.** In our theoretical analysis, we assume the following regularity conditions throughout:

(A1) We use $\lambda_{\min}(\mathbf{M})$ and $\lambda_{\max}(\mathbf{M})$ to denote the minimum and maximum eigenvalues of a positive definite matrix $\mathbf{M}$, respectively. Then, we assume

$$b \leq \lambda_{\min}\left(\frac{1}{n}\mathbf{X}^T\mathbf{X}\right) \leq \lambda_{\max}\left(\frac{1}{n}\mathbf{X}^T\mathbf{X}\right) \leq B,$$

where $b$ and $B$ are two positive constants.

(A2) $\lim_{n \to \infty} \frac{\max_{i=1,2,\ldots,n} \sum_{j=1}^p x_{ij}^2}{n} = 0$;

(A3) $E[|\varepsilon|^{2+\delta}] < \infty$ for some $\delta > 0$;

(A4) $\lim_{n \to \infty} \frac{\log(p)}{\log(n)} = \nu$ for some $0 \leq \nu < 1$.

To construct the adaptive weights ($\hat{\omega}$), we take a fixed $\gamma$ such that $\gamma > \frac{2\nu}{1-\nu}$. In our numerical studies, we let $\gamma = \lceil \frac{2\nu}{1-\nu} \rceil + 1$ to avoid the tuning on $\gamma$. Once $\gamma$ is chosen, we choose the regularization parameters according to the following conditions:

(A5) $$\lim_{n\to\infty} \frac{\lambda_2}{n} = 0, \qquad \lim_{n\to\infty} \frac{\lambda_1}{\sqrt{n}} = 0$$

and

$$\lim_{n\to\infty} \frac{\lambda_1^*}{\sqrt{n}} = 0, \qquad \lim_{n\to\infty} \frac{\lambda_1^*}{\sqrt{n}} n^{((1-\nu)(1+\gamma)-1)/2} = \infty.$$

(A6) $$\lim_{n\to\infty} \frac{\lambda_2}{\sqrt{n}} \sqrt{\sum_{j\in\mathcal{A}} \beta_j^{*2}} = 0,$$

$$\lim_{n\to\infty} \min\left(\frac{n}{\lambda_1\sqrt{p}}, \left(\frac{\sqrt{n}}{\sqrt{p}\lambda_1^*}\right)^{1/\gamma}\right)\left(\min_{j\in\mathcal{A}}|\beta_j^*|\right) \to \infty.$$

Conditions (A1) and (A2) assume the predictor matrix has a reasonably good behavior. Similar conditions were considered in Portnoy (1984). Note that in the linear regression setting, condition (A1) is exactly condition (F) in Fan and Peng (2004). Condition (A3) is used to establish the asymptotic normality of $\widehat{\boldsymbol{\beta}}$(AdaEnet).

It is worth pointing out that condition (A4) is weaker than that used in Fan and Peng (2004), in which $p$ is assumed to satisfy $p^4/n \to 0$ or at most $p^3/n \to 0$. It means their results require $\nu < \frac{1}{3}$. Our theory removes this limitation. For any $0 \le \nu < 1$, we can choose an appropriate $\gamma$ to construct the adaptive weights and the oracle property holds as long as $\gamma > \frac{2\nu}{1-\nu}$. Also note that, in the finite dimension setting, $\nu = 0$; thus, any positive $\gamma$ can be used, which agrees with the results in Zou (2006).

Condition (A6) is similar to condition (H) in Fan and Peng (2004). Basically, condition (A6) allows the nonzero coefficients to vanish but at a rate that can be distinguished by the penalized least squares. In the finite dimension setting, the condition is implicitly assumed.

THEOREM 3.1. *Given the data* $(\mathbf{y}, \mathbf{X})$, *let* $\hat{\mathbf{w}} = (\hat{w}_1, \ldots, \hat{w}_p)$ *be a vector whose components are all nonnegative and can depend on* $(\mathbf{y}, \mathbf{X})$. *Define*

$$\widehat{\boldsymbol{\beta}}_{\hat{\mathbf{w}}}(\lambda_2, \lambda_1) = \arg\min_{\boldsymbol{\beta}}\left\{\|\mathbf{y} - \mathbf{X}\boldsymbol{\beta}\|_2^2 + \lambda_2\|\boldsymbol{\beta}\|_2^2 + \lambda_1 \sum_{j=1}^p \hat{w}_j|\beta_j|\right\}$$

*for nonnegative parameters* $\lambda_2$ *and* $\lambda_1$. *If* $\hat{w}_j = 1$ *for all* $j$, *we denote* $\widehat{\boldsymbol{\beta}}_{\hat{\mathbf{w}}}(\lambda_2, \lambda_1)$ *by* $\widehat{\boldsymbol{\beta}}(\lambda_2, \lambda_1)$ *for convenience.*

*If we assume the model* (1.1) *and condition* (A1), *then*

$$E(\|\widehat{\boldsymbol{\beta}}_{\hat{\mathbf{w}}}(\lambda_2, \lambda_1) - \boldsymbol{\beta}^*\|_2^2) \le 4 \frac{\lambda_2^2\|\boldsymbol{\beta}^*\|_2^2 + Bpn\sigma^2 + \lambda_1^2 E(\sum_{j=1}^p \hat{w}_j^2)}{(bn + \lambda_2)^2}.$$

*In particular, when* $\hat{w}_j = 1$ *for all* $j$, *we have*

$$E(\|\widehat{\boldsymbol{\beta}}(\lambda_2, \lambda_1) - \boldsymbol{\beta}^*\|_2^2) \le 4 \frac{\lambda_2^2\|\boldsymbol{\beta}^*\|_2^2 + Bpn\sigma^2 + \lambda_1^2 p}{(bn + \lambda_2)^2}.$$

It is worth mentioning that the derived risk bounds are nonasymptotic. Theorem 3.1 is very useful for the asymptotic analysis. A direct corollary of Theorem 3.1 is that, under conditions (A1)–(A6), $\widehat{\boldsymbol{\beta}}(\lambda_2, \lambda_1)$ is a root-$(n/p)$-consistent estimator. This consistent rate is the same as the result of SCAD [Fan and Peng (2004)]. The root-$(n/p)$ consistency result suggests that it is appropriate to use the elastic-net to construct the adaptive weights.

THEOREM 3.2. *Let us write $\boldsymbol{\beta}^* = (\boldsymbol{\beta}^*_{\mathcal{A}}, 0)$ and define*

$$(3.1) \qquad \widetilde{\boldsymbol{\beta}}^*_{\mathcal{A}} = \arg\min_{\boldsymbol{\beta}} \left\{ \|\mathbf{y} - \mathbf{X}_{\mathcal{A}}\boldsymbol{\beta}\|_2^2 + \lambda_2 \sum_{j \in \mathcal{A}} \beta_j^2 + \lambda_1^* \sum_{j \in \mathcal{A}} \hat{w}_j |\beta_j| \right\}.$$

*Then, with probability tending to 1, $((1 + \frac{\lambda_2}{n})\widetilde{\boldsymbol{\beta}}^*_{\mathcal{A}}, 0)$ is the solution to (2.2).*

Theorem 3.2 provides an asymptotic characterization of the solution to the adaptive elastic-net criterion. The definition of $\widetilde{\boldsymbol{\beta}}^*_{\mathcal{A}}$ borrows the concept of "oracle" [Donoho and Johnstone (1994), Fan and Li (2001), Fan and Peng (2004) and Zou (2006)]. If there was an oracle informing us the true subset model, then we would use this oracle information and the adaptive elastic-net criterion would become that in (2.3). Theorem 3.2 tells us that, asymptotically speaking, the adaptive elastic-net works as if it had such oracle information. Theorem 3.2 also suggests that the adaptive elastic-net should enjoy the oracle property, which is confirmed in the next theorem.

THEOREM 3.3. *Under conditions* (A1)–(A6), *the adaptive elastic-net has the oracle property; that is, the estimator $\widehat{\boldsymbol{\beta}}(\mathrm{AdaEnet})$ must satisfy:*

1. *Consistency in selection:* $\Pr(\{j : \hat{\beta}(\mathrm{AdaEnet})_j \neq 0\} = \mathcal{A}) \to 1$,
2. *Asymptotic normality:* $\boldsymbol{\alpha}^T \frac{\mathbf{I} + \lambda_2 \boldsymbol{\Sigma}_{\mathcal{A}}^{-1}}{1 + \lambda_2/n} \boldsymbol{\Sigma}_{\mathcal{A}}^{1/2} (\widehat{\boldsymbol{\beta}}(\mathrm{AdaEnet})_{\mathcal{A}} - \boldsymbol{\beta}^*_{\mathcal{A}}) \to_d N(0, \sigma^2)$,
   *where $\boldsymbol{\Sigma}_{\mathcal{A}} = \mathbf{X}_{\mathcal{A}}^T \mathbf{X}_{\mathcal{A}}$ and $\boldsymbol{\alpha}$ is a vector of norm 1.*

By Theorem 3.3, the selection consistency and the asymptotic normality of the adaptive elastic-net are still valid when the number of parameters diverges. Technically speaking, the selection consistency result is stronger than that Theorem 3.2 implies, although Theorem 3.2 plays an important role in the proof of Theorem 3.3. As a special case, when we let $\lambda_2 = 0$, which is a choice satisfying conditions (A5) and (A6), Theorem 3.3 tells us that the adaptive lasso enjoys the selection consistency and the asymptotical normality

$$\boldsymbol{\alpha}^T \boldsymbol{\Sigma}_{\mathcal{A}}^{1/2} (\widehat{\boldsymbol{\beta}}(\mathrm{AdaLasso})_{\mathcal{A}} - \boldsymbol{\beta}^*_{\mathcal{A}}) \xrightarrow{d} N(0, \sigma^2).$$

**4. Numerical studies.** In this section, we present simulations to study the finite sample performance of the adaptive elastic-net. We considered five methods in the simulation study: the lasso (Lasso), the elastic-net (Enet), the adaptive lasso (ALasso), the adaptive elastic-net (AEnet) and the SCAD. In our implementation, we let $\lambda_2 = 0$ in the adaptive elastic-net to get the adaptive lasso fit. There are several commonly used tuning parameter selection methods, such as cross-validation, generalized cross-validation (GCV), AIC and BIC. Zou, Hastie and Tibshirani (2007) suggested using BIC to select the lasso tuning parameter. Wang, Li and Tsai (2007) showed that for the SCAD, BIC is a better tuning parameter selector than GCV and AIC. In this work, we used BIC to select the tuning parameter for each method.

Fan and Peng (2004) considered simulation models in which $p_n = [4n^{1/4}] - 5$ and $|\mathcal{A}| = 5$. Our theory allows $p_n = O(n^\nu)$ for any $\nu < 1$. Thus, we are interested in models in which $p_n = O(n^\nu)$ with $\nu > \frac{1}{3}$. In addition, we allow the intrinsic dimension ($\mathcal{A}$) to diverge with the sample size as well, because such designs make the model selection and estimation more challenging than in the fixed $|\mathcal{A}|$ situations.

EXAMPLE 1. We generated data from the linear regression model

$$y = \mathbf{x}^T \boldsymbol{\beta}^* + \varepsilon,$$

where $\boldsymbol{\beta}^*$ is a $p$-dim vector and $\varepsilon \sim N(0, \sigma^2)$, $\sigma = 6$, and $\mathbf{x}$ follows a $p$-dim multivariate normal distribution with zero mean and covariance $\boldsymbol{\Sigma}$ whose $(j, k)$ entry is $\boldsymbol{\Sigma}_{j,k} = \rho^{|j-k|}$, $1 \le k, j \le p$. We considered $\rho = 0.5$ and $\rho = 0.75$. Let $p = p_n = [4n^{1/2}] - 5$ for $n = 100, 200, 400$. Let $\mathbf{1}_m/\mathbf{0}_m$ denote a $m$-vector of 1's/0's. The true coefficients are $\boldsymbol{\beta}^* = (3 \cdot \mathbf{1}_q, 3 \cdot \mathbf{1}_q, 3 \cdot \mathbf{1}_q, \mathbf{0}_{p-3q})^T$ and $|\mathcal{A}| = 3q$ and $q = [p_n/9]$. In this example $\nu = \frac{1}{2}$; hence, we used $\gamma = 3$ for computing the adaptive weights in the adaptive elastic-net.

For each estimator $\widehat{\boldsymbol{\beta}}$, its estimation accuracy is measured by the mean squared error (MSE) defined as $E[(\widehat{\boldsymbol{\beta}} - \boldsymbol{\beta}^*)^T \boldsymbol{\Sigma} (\widehat{\boldsymbol{\beta}} - \boldsymbol{\beta}^*)]$. The variable selection performance is gauged by $(C, IC)$, where $C$ is the number of zero coefficients that are correctly estimated by zero and $IC$ is the number of nonzero coefficients that are incorrectly estimated by zero.

Table 1 documents the simulation results. Several interesting observations can be made:

1. When the sample size is large ($n = 400$), the three oracle-like estimators outperform the lasso and the elastic-net which do not have the oracle property. That is expected according to the asymptotic theory.
2. The SCAD and the adaptive elastic-net are the best when the sample size is large and the correlation is moderate. However, the SCAD can perform much worse than the adaptive elastic-net when the correlation is high ($\rho = 0.75$) or the sample size is small.

3. Both the elastic-net and the adaptive lasso can do significantly better than the lasso. What is more interesting is that the adaptive elastic-net often outperforms the elastic-net and the adaptive lasso.

EXAMPLE 2. We considered the same setup as in Example 1, except that we let $p = p_n = [4n^{2/3}] - 5$ for $n = 100, 200, 800$. Since $\nu = \frac{2}{3}$, we used $\gamma = 5$ for computing the adaptive weights in the adaptive elastic-net and the adaptive lasso. The estimation problem in this example is even more difficult than that in Example 1. To see why, note that when $n = 200$ the dimension increases from 51 in Example 1 to 131 in this example, and the intrinsic dimension ($|\mathcal{A}|$) is almost tripled.

The simulation results are presented in Table 2, from which we can see that the three observations made in Example 1 are still valid in this example. Furthermore, we see that, for every combination of $(n, p, |\mathcal{A}|, \rho)$, the adaptive elastic-net has the best performance.

**5. Ultra-high dimensional data.** In this section, we discuss how the adaptive elastic-net can be applied to ultra-high dimensional data in which $p > n$. When $p$ is much larger than $n$, Candes and Tao (2007) suggested using the Dantzig selector which can achieve the ideal estimation risk up to a $\log(p)$ factor under the uniform uncertainty condition. Fan and Lv (2008) showed that the uniform uncertainty condition may easily fail and the $\log(p)$ factor is too large when $p$ is exponentially large. Moreover, the computational cost of the Dantzig selector would be very high when $p$ is large. In order to overcome these difficulties, Fan and Lv (2008) introduced the Sure Independence Screening (SIS) idea, which reduces the ultra-high dimensionality to a relatively large scale $d_n$ but $d_n < n$. Then, the lower dimension methods such as the SCAD can be used to estimate the sparse model. This procedure is referred to as SIS + SCAD. Under regularity conditions, Fan and Lv (2008) proved that SIS misses true features with an exponentially small probability and SIS + SCAD holds the oracle property if $d_n = o(n^{1/3})$. Furthermore, with the help of SIS, the Dantzig selector can achieve the ideal risk up to a $\log(d_n)$ factor, rather than the original $\log(p)$.

Inspired by the results of Fan and Lv (2008), we consider combining the adaptive elastic-net and SIS when $p > n$. We first apply SIS to reduce the dimension to $d_n$ and then fit the data by using the adaptive elastic-net. We call this procedure SIS + AEnet.

THEOREM 5.1. *Suppose the conditions for Theorem 1 in Fan and Lv (2008) hold. Let $d_n = O(n^\nu)$, $\nu < 1$; then, SIS+AEnet produces an estimator that holds the oracle property.*

Table 1
*Simulation* I: *model selection and fitting results based on 100 replications*

| $n$ | $p_n$ | $|\mathcal{A}|$ | Model | MSE | C | IC |
|---|---|---|---|---|---|---|
| | | | $\rho = 0.5$ | | | |
| 100 | 35 | 9 | Truth | | 26 | 0 |
| | | | Lasso | 7.57 (0.31) | 24.08 | 0.01 |
| | | | ALasso | 6.78 (0.42) | 25.50 | 0.42 |
| | | | Enet | 5.91 (0.29) | 24.06 | 0 |
| | | | AEnet | 5.07 (0.35) | 25.47 | 0.15 |
| | | | SCAD | 10.55 (0.68) | 22.54 | 0.35 |
| 200 | 51 | 15 | Truth | | 36 | 0 |
| | | | Lasso | 6.63 (0.24) | 33.32 | 0 |
| | | | ALasso | 3.78 (0.18) | 35.46 | 0.02 |
| | | | Enet | 4.86 (0.19) | 33.36 | 0 |
| | | | AEnet | 3.46 (0.17) | 35.47 | 0.01 |
| | | | SCAD | 4.76 (0.33) | 34.63 | 0.10 |
| 400 | 75 | 24 | Truth | | 51 | 0 |
| | | | Lasso | 4.99 (0.15) | 47.31 | 0 |
| | | | ALasso | 2.76 (0.09) | 50.33 | 0 |
| | | | Enet | 3.37 (0.12) | 48.00 | 0 |
| | | | AEnet | 2.47 (0.08) | 50.45 | 0 |
| | | | SCAD | 2.42 (0.09) | 50.88 | 0 |
| | | | $\rho = 0.75$ | | | |
| 100 | 35 | 9 | Truth | | 26 | 0 |
| | | | Lasso | 5.93 (0.26) | 24.80 | 0.14 |
| | | | ALasso | 8.49 (0.39) | 25.76 | 1.84 |
| | | | Enet | 4.18 (0.24) | 24.77 | 0.05 |
| | | | AEnet | 5.24 (0.32) | 25.70 | 0.74 |
| | | | SCAD | 11.59 (0.56) | 22.46 | 1.34 |
| 200 | 51 | 15 | Truth | | 36 | 0 |
| | | | Lasso | 5.10 (0.18) | 34.66 | 0.02 |
| | | | ALasso | 5.32 (0.31) | 35.70 | 0.87 |
| | | | Enet | 3.79 (0.17) | 34.79 | 0 |
| | | | AEnet | 3.32 (0.17) | 35.80 | 0.19 |
| | | | SCAD | 5.99 (0.31) | 33.10 | 0.35 |
| 400 | 75 | 24 | Truth | | 51 | 0 |
| | | | Lasso | 3.83 (0.12) | 49.03 | 0 |
| | | | ALasso | 2.85 (0.12) | 50.53 | 0.09 |
| | | | Enet | 3.24 (0.11) | 49.07 | 0 |
| | | | AEnet | 2.71 (0.09) | 50.54 | 0.03 |
| | | | SCAD | 3.64 (0.17) | 48.43 | 0.09 |

We make a note here that Theorem 5.1 is a direct consequence of Theorem 1 in Fan and Lv (2008) and Theorem 3.3; thus, its proof is omitted. Theorem 5.1 is similar to Theorem 5 in Fan and Lv (2008), but there is a difference. SIS + AEnent can hold the oracle property when $d_n$ exceeds $O(n^{1/3})$, while Theorem 5 in Fan and Lv (2008) assumes $d_n = o(n^{1/3})$.

Table 2
*Example 2: model selection and fitting results based on 100 replications*

| $n$ | $p_n$ | $|\mathcal{A}|$ | Model | MSE | C | IC |
|---|---|---|---|---|---|---|
| | | | $\rho = 0.5$ | | | |
| 100 | 81 | 27 | Truth | | 54 | 0 |
| | | | Lasso | 31.73 (1.06) | 47.06 | 0.19 |
| | | | ALasso | 28.78 (1.22) | 53.01 | 2.12 |
| | | | Enet | 27.61 (1.04) | 46.35 | 0.13 |
| | | | AEnet | 20.27 (0.94) | 53.00 | 1.15 |
| | | | SCAD | 44.88 (2.65) | 47.79 | 2.37 |
| 200 | 131 | 42 | Truth | | 89 | 0 |
| | | | Lasso | 23.41 (0.67) | 80.51 | 0 |
| | | | ALasso | 12.70 (0.48) | 87.99 | 0.14 |
| | | | Enet | 18.94 (0.61) | 80.27 | 0 |
| | | | AEnet | 10.68 (0.37) | 87.97 | 0 |
| | | | SCAD | 14.14 (0.64) | 87.42 | 0.25 |
| 800 | 339 | 111 | Truth | | 228 | 0 |
| | | | Lasso | 13.72 (0.23) | 212.10 | 0 |
| | | | ALasso | 6.44 (0.12) | 226.61 | 0 |
| | | | Enet | 11.02 (0.18) | 213.91 | 0 |
| | | | AEnet | 6.00 (0.10) | 226.75 | 0 |
| | | | SCAD | 7.79 (0.30) | 228.00 | 0.33 |
| | | | $\rho = 0.75$ | | | |
| 100 | 81 | 27 | Truth | | 54 | 0 |
| | | | Lasso | 22.04 (0.73) | 50.74 | 0.71 |
| | | | ALasso | 33.98 (1.08) | 53.73 | 7.19 |
| | | | Enet | 17.37 (0.62) | 50.82 | 0.46 |
| | | | AEnet | 16.18 (0.80) | 53.67 | 2.36 |
| | | | SCAD | 31.84 (1.77) | 50.55 | 4.74 |
| 200 | 131 | 42 | Truth | | 89 | 0 |
| | | | Lasso | 16.71 (0.50) | 85.17 | 0.06 |
| | | | ALasso | 20.98 (0.92) | 88.64 | 3.98 |
| | | | Enet | 14.12 (0.48) | 85.35 | 0.05 |
| | | | AEnet | 11.16 (0.46) | 88.60 | 0.87 |
| | | | SCAD | 15.27 (0.61) | 87.20 | 1.33 |
| 800 | 339 | 111 | Truth | | 228 | 0 |
| | | | Lasso | 10.01 (0.16) | 221.74 | 0 |
| | | | ALasso | 6.39 (0.12) | 226.89 | 0 |
| | | | Enet | 8.01 (0.13) | 222.74 | 0 |
| | | | AEnet | 6.23 (0.11) | 226.94 | 0 |
| | | | SCAD | 6.62 (0.17) | 228.00 | 0.29 |

TABLE 3
*A demonstration of SIS + AEnet: model selection and fitting results based on 100 replications*

| $d_n = [5.5n^{2/3}]$ | Model | MSE | C | IC |
|---|---|---|---|---|
| 188 | Truth | | 992 | 0 |
| | SIS + AEnet | 0.71 (0.18) | 987.45 | 0.05 |
| | SIS + SCAD | 1.48 (0.90) | 982.20 | 0.06 |

To demonstrate SIS + AEnet, we consider the simulation example used in Fan and Lv (2008), Section 3.3.1. The model is $y = \mathbf{x}^T \boldsymbol{\beta}^* + 1.5 N(0,1)$, where $\boldsymbol{\beta}^* = (\boldsymbol{\beta}_1^T, \mathbf{0}_{p-|\mathcal{A}|})^T$ with $|\mathcal{A}| = 8$. Here, $\boldsymbol{\beta}_1$ is a 8-dim vector and each component has the form $(-1)^u (a_n + |z|)$, where $a_n = 4\log(n)/\sqrt{n}$, $u$ is randomly drawn from Ber(0.4) and $z$ is randomly drawn from the standard normal distribution. We generated $n = 200$ data from the above model. Before applying the adaptive elastic-net, we used SIS to reduce the dimensionality from 1000 to $d_n = [5.5n^{2/3}] = 188$. The estimation problem is still rather challenging, as we need to estimate 188 parameters by using only 200 observations. From Table 3, we see that SIS + AEnet performs favorably compared to SIS + SCAD.

## 6. Proofs.

PROOF OF THEOREM 3.1. We write
$$\widehat{\boldsymbol{\beta}}(\lambda_2, 0) = \arg\min_{\boldsymbol{\beta}} \|\mathbf{y} - \mathbf{X}\boldsymbol{\beta}\|_2^2 + \lambda_2 \|\boldsymbol{\beta}\|_2^2.$$

By the definition of $\widehat{\boldsymbol{\beta}}_{\widehat{\mathbf{w}}}(\lambda_2, \lambda_1)$ and $\widehat{\boldsymbol{\beta}}(\lambda_2, 0)$, we know

$$\|\mathbf{y} - \mathbf{X}\widehat{\boldsymbol{\beta}}_{\widehat{\mathbf{w}}}(\lambda_2, \lambda_1)\|_2^2 + \lambda_2 \|\widehat{\boldsymbol{\beta}}_{\widehat{\mathbf{w}}}(\lambda_2, \lambda_1)\|_2^2 \geq \|\mathbf{y} - \mathbf{X}\widehat{\boldsymbol{\beta}}(\lambda_2, 0)\|_2^2 + \lambda_2 \|\widehat{\boldsymbol{\beta}}(\lambda_2, 0)\|_2^2$$

and

$$\|\mathbf{y} - \mathbf{X}\widehat{\boldsymbol{\beta}}(\lambda_2, 0)\|_2^2 + \lambda_2 \|\widehat{\boldsymbol{\beta}}(\lambda_2, 0)\|_2^2 + \lambda_1 \sum_{j=1}^{p} \widehat{w}_j |\hat{\beta}(\lambda_2, 0)_j|$$
$$\geq \|\mathbf{y} - \mathbf{X}\widehat{\boldsymbol{\beta}}_{\widehat{\mathbf{w}}}(\lambda_2, \lambda_1)\|_2^2 + \lambda_2 \|\widehat{\boldsymbol{\beta}}_{\widehat{\mathbf{w}}}(\lambda_2, \lambda_1)\|_2^2 + \lambda_1 \sum_{j=1}^{p} \widehat{w}_j |\hat{\beta}_{\widehat{\mathbf{w}}}(\lambda_2, \lambda_1)_j|.$$

From the above two inequalities, we have

$$\lambda_1 \sum_{j=1}^{p} \widehat{w}_j (|\hat{\beta}(\lambda_2, 0)_j| - |\hat{\beta}_{\widehat{\mathbf{w}}}(\lambda_2, \lambda_1)_j|)$$
(6.1)
$$\geq (\|\mathbf{y} - \mathbf{X}\widehat{\boldsymbol{\beta}}_{\widehat{\mathbf{w}}}(\lambda_2, \lambda_1)\|_2^2 + \lambda_2 \|\widehat{\boldsymbol{\beta}}_{\widehat{\mathbf{w}}}(\lambda_2, \lambda_1)\|_2^2)$$
$$- (\|\mathbf{y} - \mathbf{X}\widehat{\boldsymbol{\beta}}(\lambda_2, 0)\|_2^2 + \lambda_2 \|\widehat{\boldsymbol{\beta}}(\lambda_2, 0)\|_2^2).$$

On the other hand, we have

$$(\|\mathbf{y} - \mathbf{X}\widehat{\boldsymbol{\beta}}_{\widehat{\mathbf{w}}}(\lambda_2, \lambda_1)\|_2^2 + \lambda_2\|\widehat{\boldsymbol{\beta}}_{\widehat{\mathbf{w}}}(\lambda_2, \lambda_1)\|_2^2)$$
$$- (\|\mathbf{y} - \mathbf{X}\widehat{\boldsymbol{\beta}}(\lambda_2, 0)\|_2^2 + \lambda_2\|\widehat{\boldsymbol{\beta}}(\lambda_2, 0)\|_2^2)$$
$$= (\widehat{\boldsymbol{\beta}}_{\widehat{\mathbf{w}}}(\lambda_2, \lambda_1) - \widehat{\boldsymbol{\beta}}(\lambda_2, 0))^T (\mathbf{X}^T\mathbf{X} + \lambda_2\mathbf{I})(\widehat{\boldsymbol{\beta}}_{\widehat{\mathbf{w}}}(\lambda_2, \lambda_1) - \widehat{\boldsymbol{\beta}}(\lambda_2, 0))$$

and

$$\sum_{j=1}^{p} \hat{w}_j (|\hat{\beta}(\lambda_2, 0)_j| - |\hat{\beta}_{\widehat{\mathbf{w}}}(\lambda_2, \lambda_1)_j|)$$
$$\leq \sum_{j=1}^{p} \hat{w}_j |\hat{\beta}(\lambda_2, 0)_j - \hat{\beta}_{\widehat{\mathbf{w}}}(\lambda_2, \lambda_1)_j|$$
$$\leq \sqrt{\sum_{j=1}^{p} \hat{w}_j^2} \|\widehat{\boldsymbol{\beta}}(\lambda_2, 0) - \widehat{\boldsymbol{\beta}}_{\widehat{\mathbf{w}}}(\lambda_2, \lambda_1)\|_2.$$

Note that $\lambda_{\min}(\mathbf{X}^T\mathbf{X} + \lambda_2\mathbf{I}) = \lambda_{\min}(\mathbf{X}^T\mathbf{X}) + \lambda_2$. Therefore, we end up with

$$(\lambda_{\min}(\mathbf{X}^T\mathbf{X}) + \lambda_2)\|\widehat{\boldsymbol{\beta}}_{\widehat{\mathbf{w}}}(\lambda_2, \lambda_1) - \widehat{\boldsymbol{\beta}}(\lambda_2, 0)\|_2^2$$
(6.2) $$\leq (\widehat{\boldsymbol{\beta}}_{\widehat{\mathbf{w}}}(\lambda_2, \lambda_1) - \widehat{\boldsymbol{\beta}}(\lambda_2, 0))^T (\mathbf{X}^T\mathbf{X} + \lambda_2\mathbf{I})(\widehat{\boldsymbol{\beta}}_{\widehat{\mathbf{w}}}(\lambda_2, \lambda_1) - \widehat{\boldsymbol{\beta}}(\lambda_2, 0))$$
$$\leq \lambda_1 \sqrt{\sum_{j=1}^{p} \hat{w}_j^2} \|\widehat{\boldsymbol{\beta}}(\lambda_2, 0) - \widehat{\boldsymbol{\beta}}_{\widehat{\mathbf{w}}}(\lambda_2, \lambda_1)\|_2,$$

which results in the inequality

(6.3) $$\|\widehat{\boldsymbol{\beta}}_{\widehat{\mathbf{w}}}(\lambda_2, \lambda_1) - \widehat{\boldsymbol{\beta}}(\lambda_2, 0)\|_2 \leq \frac{\lambda_1 \sqrt{\sum_{j=1}^{p} \hat{w}_j^2}}{\lambda_{\min}(\mathbf{X}^T\mathbf{X}) + \lambda_2}.$$

Note that

$$\widehat{\boldsymbol{\beta}}(\lambda_2, 0) - \boldsymbol{\beta}^* = -\lambda_2(\mathbf{X}^T\mathbf{X} + \lambda_2\mathbf{I})^{-1}\boldsymbol{\beta}^* + (\mathbf{X}^T\mathbf{X} + \lambda_2\mathbf{I})^{-1}\mathbf{X}^T\boldsymbol{\varepsilon},$$

which implies that

$$E(\|\widehat{\boldsymbol{\beta}}(\lambda_2, 0) - \boldsymbol{\beta}^*\|_2^2)$$
$$\leq 2\lambda_2^2 \|(\mathbf{X}^T\mathbf{X} + \lambda_2\mathbf{I})^{-1}\boldsymbol{\beta}^*\|_2^2 + 2E(\|(\mathbf{X}^T\mathbf{X} + \lambda_2\mathbf{I})^{-1}\mathbf{X}^T\boldsymbol{\varepsilon}\|_2^2)$$
$$\leq 2\lambda_2^2 (\lambda_{\min}(\mathbf{X}^T\mathbf{X}) + \lambda_2)^{-2} \|\boldsymbol{\beta}^*\|_2^2$$
(6.4)
$$\quad + 2(\lambda_{\min}(\mathbf{X}^T\mathbf{X}) + \lambda_2)^{-2} E(\boldsymbol{\varepsilon}^T\mathbf{X}\mathbf{X}^T\boldsymbol{\varepsilon})$$
$$= 2(\lambda_{\min}(\mathbf{X}^T\mathbf{X}) + \lambda_2)^{-2} (\lambda_2^2 \|\boldsymbol{\beta}^*\|_2^2 + \text{Tr}(\mathbf{X}^T\mathbf{X})\sigma^2)$$
$$\leq 2(\lambda_{\min}(\mathbf{X}^T\mathbf{X}) + \lambda_2)^{-2} (\lambda_2^2 \|\boldsymbol{\beta}^*\|_2^2 + p\lambda_{\max}(\mathbf{X}^T\mathbf{X})\sigma^2).$$

Combing (6.3) and (6.4), we have

$$E(\|\widehat{\boldsymbol{\beta}}_{\widehat{\mathbf{w}}}(\lambda_2, \lambda_1) - \boldsymbol{\beta}^*\|_2^2)$$
$$\leq 2E(\|\widehat{\boldsymbol{\beta}}(\lambda_2, 0) - \boldsymbol{\beta}^*\|_2^2) + 2E(\|\widehat{\boldsymbol{\beta}}_{\widehat{\mathbf{w}}}(\lambda_2, \lambda_1) - \widehat{\boldsymbol{\beta}}(\lambda_2, 0)\|_2^2)$$

$$(6.5) \qquad \leq \frac{4\lambda_2^2\|\boldsymbol{\beta}^*\|_2^2 + 4p\lambda_{\max}(\mathbf{X}^T\mathbf{X})\sigma^2 + 2\lambda_1^2 E[\sum_{j=1}^p \hat{w}_j^2]}{(\lambda_{\min}(\mathbf{X}^T\mathbf{X}) + \lambda_2)^2}$$

$$(6.6) \qquad \leq 4\frac{\lambda_2^2\|\boldsymbol{\beta}^*\|_2^2 + Bpn\sigma^2 + \lambda_1^2 E[\sum_{j=1}^p \hat{w}_j^2]}{(bn + \lambda_2)^2}.$$

We have used condition (A1) in the last inequality. When $\hat{w}_j = 1$ for all $j$, we have

$$E(\|\widehat{\boldsymbol{\beta}}(\lambda_2, \lambda_1) - \boldsymbol{\beta}^*\|_2^2) \leq 4\frac{\lambda_2^2\|\boldsymbol{\beta}^*\|_2^2 + Bpn\sigma^2 + p\lambda_1^2}{(bn + \lambda_2)^2}. \qquad \square$$

PROOF OF THEOREM 3.2. We show that $((1 + \frac{\lambda_2}{n})\widetilde{\boldsymbol{\beta}}_{\mathcal{A}}^*, 0)$ satisfies the Karush–Kuhn–Tucker (KKT) conditions of (2.2) with probability tending to 1. By the definition of $\widetilde{\boldsymbol{\beta}}_{\mathcal{A}}^*$, it suffices to show

$$\Pr(\forall j \in \mathcal{A}^c \ |-2X_j^T(\mathbf{y} - \mathbf{X}_{\mathcal{A}}\widetilde{\boldsymbol{\beta}}_{\mathcal{A}}^*)| \leq \lambda_1^* \hat{w}_j) \to 1$$

or, equivalently,

$$\Pr(\exists j \in \mathcal{A}^c \ |-2X_j^T(\mathbf{y} - \mathbf{X}_{\mathcal{A}}\widetilde{\boldsymbol{\beta}}_{\mathcal{A}}^*)| > \lambda_1^* \hat{w}_j) \to 0.$$

Let $\eta = \min_{j \in \mathcal{A}}(|\beta_j^*|)$ and $\hat{\eta} = \min_{j \in \mathcal{A}}(|\hat{\beta}(\text{enet})_j^*|)$. We note that

$$\Pr(\exists j \in \mathcal{A}^c \ |-2X_j^T(\mathbf{y} - \mathbf{X}_{\mathcal{A}}\widetilde{\boldsymbol{\beta}}_{\mathcal{A}}^*)| > \lambda_1^* \hat{w}_j)$$
$$\leq \sum_{j \in \mathcal{A}^c} \Pr(|-2X_j^T(\mathbf{y} - \mathbf{X}_{\mathcal{A}}\widetilde{\boldsymbol{\beta}}_{\mathcal{A}}^*)| > \lambda_1^* \hat{w}_j, \hat{\eta} > \eta/2) + \Pr(\hat{\eta} \leq \eta/2),$$

$$\Pr(\hat{\eta} \leq \eta/2) \leq \Pr(\|\widehat{\boldsymbol{\beta}}(\text{enet}) - \boldsymbol{\beta}^*\|_2 \geq \eta/2) \leq \frac{E(\|\widehat{\boldsymbol{\beta}}(\text{enet}) - \boldsymbol{\beta}^*\|_2^2)}{\eta^2/4}.$$

Then, by Theorem 3.1, we obtain

$$(6.7) \qquad \Pr(\hat{\eta} \leq \eta/2) \leq 16\frac{\lambda_2^2\|\boldsymbol{\beta}^*\|_2^2 + Bpn\sigma^2 + \lambda_1^2 p}{(bn + \lambda_2)^2 \eta^2}.$$

Moreover, let $M = (\frac{\lambda_1^*}{n})^{1/(1+\gamma)}$, and we have

$$\sum_{j \in \mathcal{A}^c} \Pr(|-2X_j^T(\mathbf{y} - \mathbf{X}_{\mathcal{A}}\widetilde{\boldsymbol{\beta}}_{\mathcal{A}}^*)| > \lambda_1^* \hat{w}_j, \hat{\eta} > \eta/2)$$
$$\leq \sum_{j \in \mathcal{A}^c} \Pr(|-2X_j^T(\mathbf{y} - \mathbf{X}_{\mathcal{A}}\widetilde{\boldsymbol{\beta}}_{\mathcal{A}}^*)| > \lambda_1^* \hat{w}_j, \hat{\eta} > \eta/2, |\hat{\beta}(\text{enet})_j| \leq M)$$

$$+ \sum_{j \in \mathcal{A}^c} \Pr(|\hat{\beta}(\text{enet})_j| > M)$$

$$\leq \sum_{j \in \mathcal{A}^c} \Pr(|-2X_j^T(\mathbf{y} - \mathbf{X}_\mathcal{A} \widetilde{\boldsymbol{\beta}}_\mathcal{A}^*)| > \lambda_1^* M^{-\gamma}, \hat{\eta} > \eta/2)$$

$$+ \sum_{j \in \mathcal{A}^c} \Pr(|\hat{\beta}(\text{enet})_j| > M)$$

$$(6.8) \quad \leq \frac{4M^{2\gamma}}{\lambda_1^{*2}} E\left( \sum_{j \in \mathcal{A}^c} |X_j^T(\mathbf{y} - \mathbf{X}_\mathcal{A} \widetilde{\boldsymbol{\beta}}_\mathcal{A}^*)|^2 I(\hat{\eta} > \eta/2) \right)$$

$$+ \frac{1}{M^2} E\left( \sum_{j \in \mathcal{A}^c} |\hat{\beta}(\text{enet})_j|^2 \right)$$

$$\leq \frac{4M^{2\gamma}}{\lambda_1^{*2}} E\left( \sum_{j \in \mathcal{A}^c} |X_j^T(\mathbf{y} - \mathbf{X}_\mathcal{A} \widetilde{\boldsymbol{\beta}}_\mathcal{A}^*)|^2 I(\hat{\eta} > \eta/2) \right)$$

$$+ \frac{E(\|\widehat{\boldsymbol{\beta}}(\text{enet}) - \boldsymbol{\beta}^*\|_2^2)}{M^2}$$

$$\leq \frac{4M^{2\gamma}}{\lambda_1^{*2}} E\left( \sum_{j \in \mathcal{A}^c} |X_j^T(\mathbf{y} - \mathbf{X}_\mathcal{A} \widetilde{\boldsymbol{\beta}}_\mathcal{A}^*)|^2 I(\hat{\eta} > \eta/2) \right)$$

$$+ 4\frac{\lambda_2^2 \|\boldsymbol{\beta}^*\|_2^2 + Bpn\sigma^2 + \lambda_1^2 p}{(bn + \lambda_2)^2 M^2},$$

where we have used Theorem 3.1 in the last step. By the model assumption, we have

$$\sum_{j \in \mathcal{A}^c} |X_j^T(\mathbf{y} - \mathbf{X}_\mathcal{A} \widetilde{\boldsymbol{\beta}}_\mathcal{A}^*)|^2 = \sum_{j \in \mathcal{A}^c} |X_j^T(\mathbf{X}_\mathcal{A} \boldsymbol{\beta}_\mathcal{A}^* - \mathbf{X}_\mathcal{A} \widetilde{\boldsymbol{\beta}}_\mathcal{A}^*) + X_j^T \boldsymbol{\varepsilon}|^2$$

$$\leq 2 \sum_{j \in \mathcal{A}^c} |X_j^T(\mathbf{X}_\mathcal{A} \boldsymbol{\beta}_\mathcal{A}^* - \mathbf{X}_\mathcal{A} \widetilde{\boldsymbol{\beta}}_\mathcal{A}^*)|^2 + 2 \sum_{j \in \mathcal{A}^c} |X_j^T \boldsymbol{\varepsilon}|^2$$

$$\leq 2Bn \|\mathbf{X}_\mathcal{A}(\boldsymbol{\beta}_\mathcal{A}^* - \widetilde{\boldsymbol{\beta}}_\mathcal{A}^*)\|_2^2 + 2 \sum_{j \in \mathcal{A}^c} |X_j^T \boldsymbol{\varepsilon}|^2$$

$$\leq 2Bn \cdot Bn \|\boldsymbol{\beta}_\mathcal{A}^* - \widetilde{\boldsymbol{\beta}}_\mathcal{A}^*\|_2^2 + 2 \sum_{j \in \mathcal{A}^c} |X_j^T \boldsymbol{\varepsilon}|^2,$$

which gives us the inequality

$$(6.9) \quad E\left( \sum_{j \in \mathcal{A}^c} |X_j^T(\mathbf{y} - \mathbf{X}_\mathcal{A} \widetilde{\boldsymbol{\beta}}_\mathcal{A}^*)|^2 I(\hat{\eta} > \eta/2) \right)$$

$$\leq 2B^2 n^2 E(\|\boldsymbol{\beta}_\mathcal{A}^* - \widetilde{\boldsymbol{\beta}}_\mathcal{A}^*\|_2^2 I(\hat{\eta} > \eta/2)) + 2Bnp\sigma^2.$$

We now bound $E(\|\boldsymbol{\beta}^*_{\mathcal{A}} - \widetilde{\boldsymbol{\beta}}^*_{\mathcal{A}}\|^2_2 I(\hat{\eta} > \eta/2))$. Let

$$\widetilde{\boldsymbol{\beta}}^*_{\mathcal{A}}(\lambda_2, 0) = \arg\min_{\boldsymbol{\beta}} \left\{ \|\mathbf{y} - \mathbf{X}_{\mathcal{A}}\boldsymbol{\beta}\|^2_2 + \lambda_2 \sum_{j \in \mathcal{A}} \beta_j^2 \right\}.$$

Then, by using the same arguments for deriving (6.1), (6.2) and (6.3), we have

$$(6.10) \quad \|\widetilde{\boldsymbol{\beta}}^*_{\mathcal{A}} - \widetilde{\boldsymbol{\beta}}^*_{\mathcal{A}}(\lambda_2, 0)\|_2 \le \frac{\lambda_1^* \cdot \max_{j \in \mathcal{A}} \hat{w}_j \sqrt{|\mathcal{A}|}}{\lambda_{\min}(\mathbf{X}_{\mathcal{A}}^T\mathbf{X}_{\mathcal{A}}) + \lambda_2} \le \frac{\lambda_1^* \hat{\eta}^{-\gamma} \sqrt{p}}{bn + \lambda_2}.$$

Note that $\lambda_{\min}(\mathbf{X}_{\mathcal{A}}^T\mathbf{X}_{\mathcal{A}}) \ge \lambda_{\min}(\mathbf{X}^T\mathbf{X}) \ge bn$ and $\lambda_{\max}(\mathbf{X}_{\mathcal{A}}^T\mathbf{X}_{\mathcal{A}}) \le \lambda_{\max}(\mathbf{X}^T\mathbf{X}) \le Bn$. Following the rest arguments in the proof of Theorem 3.1, we obtain

$$(6.11) \quad \begin{aligned} E(\|\boldsymbol{\beta}^*_{\mathcal{A}} &- \widetilde{\boldsymbol{\beta}}^*_{\mathcal{A}}\|^2_2 I(\hat{\eta} > \eta/2)) \\ &\le 4\frac{\lambda_2^2\|\boldsymbol{\beta}^*_{\mathcal{A}}\|^2_2 + \lambda_{\max}(\mathbf{X}_{\mathcal{A}}^T\mathbf{X}_{\mathcal{A}})|\mathcal{A}|\sigma^2 + \lambda_1^{*2}(\eta/2)^{-2\gamma}|\mathcal{A}|}{(\lambda_{\min}(\mathbf{X}_{\mathcal{A}}^T\mathbf{X}_{\mathcal{A}}) + \lambda_2)^2} \\ &\le 4\frac{\lambda_2^2\|\boldsymbol{\beta}^*\|^2_2 + Bpn\sigma^2 + \lambda_1^{*2}(\eta/2)^{-2\gamma}p}{(bn + \lambda_2)^2}. \end{aligned}$$

The combination of (6.7), (6.8), (6.9) and (6.11) yields

$$\begin{aligned} \Pr(\exists j \in \mathcal{A}^c \ &|-2X_j^T(\mathbf{y} - \mathbf{X}_{\mathcal{A}}\widetilde{\boldsymbol{\beta}}^*_{\mathcal{A}})| > \lambda_1^* \hat{w}_j) \\ &\le \frac{4M^{2\gamma}n}{\lambda_1^{*2}} \left( 8B^2 n \frac{\lambda_2^2\|\boldsymbol{\beta}^*\|^2_2 + Bpn\sigma^2 + \lambda_1^{*2}(\eta/2)^{-2\gamma}p}{(bn + \lambda_2)^2} + 2Bp\sigma^2 \right) \\ &\quad + \frac{\lambda_2^2\|\boldsymbol{\beta}^*\|^2_2 + Bpn\sigma^2 + \lambda_1^2 p}{(bn + \lambda_2)^2} \frac{4}{M^2} + \frac{\lambda_2^2\|\boldsymbol{\beta}^*\|^2_2 + Bpn\sigma^2 + \lambda_1^2 p}{(bn + \lambda_2)^2} \frac{16}{\eta^2} \\ &\triangleq K_1 + K_2 + K_3. \end{aligned}$$

We have chosen $\gamma > \frac{2\nu}{1-\nu}$; then, under conditions (A1)–(A6), it follows that

$$\begin{aligned} K_1 &= O\left( \left( \frac{\lambda_1^*}{\sqrt{n}} n^{((1+\gamma)(1-\nu)-1)/2} \right)^{-2/(1+\gamma)} \right) \to 0, \\ K_2 &= O\left( \frac{p}{n} \left( \frac{n}{\lambda_1^*} \right)^{2/(1+\gamma)} \right) \to 0, \\ (6.12) \quad K_3 &= O\left( \frac{p}{n} \frac{1}{\eta^2} \right) \\ &= O\left( \left( \lambda_1^* \sqrt{\frac{p}{n}} \eta^{-\gamma} \right)^{2/\gamma} \left( \frac{p}{n} \left( \frac{n}{\lambda_1^*} \right)^{2/(1+\gamma)} \right)^{(1+\gamma)/\gamma} p^{-2/\gamma} \right) \to 0. \end{aligned}$$

Thus, the proof is complete. □

PROOF OF THEOREM 3.3. From Theorem 3.2, we have shown that, with probability tending to 1, the adaptive elastic-net estimator is equal to $((1+\frac{\lambda_2}{n})\widetilde{\boldsymbol{\beta}}^*_{\mathcal{A}}, 0)$. Therefore, in order to prove the model selection consistency result, we only need to show $\Pr(\min_{j \in \mathcal{A}} |\tilde{\beta}^*_j| > 0) \to 1$. By (6.10), we have

$$\min_{j \in \mathcal{A}} |\tilde{\beta}^*_j| > \min_{j \in \mathcal{A}} |\tilde{\beta}^*(\lambda_2, 0)_j| - \frac{\lambda^*_1 \sqrt{p} \hat{\eta}^{-\gamma}}{bn + \lambda_2}.$$

Note that

$$\min_{j \in \mathcal{A}} |\tilde{\beta}^*(\lambda_2, 0)_j| > \min_{j \in \mathcal{A}} |\beta^*_j| - \|\widetilde{\boldsymbol{\beta}}^*_{\mathcal{A}}(\lambda_2, 0) - \boldsymbol{\beta}^*_{\mathcal{A}}\|_2.$$

Following (6.6), it is easy to see that

$$E(\|\widetilde{\boldsymbol{\beta}}^*_{\mathcal{A}}(\lambda_2, 0) - \boldsymbol{\beta}^*_{\mathcal{A}}\|_2^2) \le 4\frac{\lambda_2^2 \|\boldsymbol{\beta}^*\|_2^2 + Bpn\sigma^2}{(bn + \lambda_2)^2} = O\left(\frac{p}{n}\right).$$

Moreover, $\frac{\lambda^*_1 \sqrt{p} \hat{\eta}^{-\gamma}}{bn+\lambda_2} = O(\frac{1}{\sqrt{n}})(\frac{\lambda^*_1 \sqrt{p}}{\sqrt{n}} \eta^{-\gamma})(\frac{\hat{\eta}}{\eta})^{-\gamma}$ and

$$E\left(\left(\frac{\hat{\eta}}{\eta}\right)^2\right) \le 2 + \frac{2}{\eta^2} E((\hat{\eta} - \eta)^2)$$

$$\le 2 + \frac{2}{\eta^2} E(\|\widehat{\boldsymbol{\beta}}(\lambda_2, \lambda_1) - \boldsymbol{\beta}^*\|_2^2)$$

$$\le 2 + \frac{8}{\eta^2} \frac{\lambda_2^2 \|\boldsymbol{\beta}^*\|_2^2 + Bpn\sigma^2 + \lambda_1^2 p}{(bn + \lambda_2)^2}.$$

In (6.12) we have shown $\eta^2 \frac{n}{p} \to \infty$. Thus,

(6.13) $$\frac{\lambda^*_1 \sqrt{p} \hat{\eta}^{-\gamma}}{bn + \lambda_2} = o\left(\frac{1}{\sqrt{n}}\right) O_P(1).$$

Hence, we have

$$\min_{j \in \mathcal{A}} |\tilde{\beta}^*_j| > \eta - \sqrt{\frac{p}{n}} O_P(1) - o\left(\frac{1}{\sqrt{n}}\right) O_P(1)$$

and $\Pr(\min_{j \in \mathcal{A}} |\tilde{\beta}^*_j| > 0) \to 1$.

We now prove the asymptotic normality. For convenience, we write

$$z_n = \boldsymbol{\alpha}^T \frac{\mathbf{I} + \lambda_2 \boldsymbol{\Sigma}^{-1}_{\mathcal{A}}}{1 + \lambda_2/n} \boldsymbol{\Sigma}^{1/2}_{\mathcal{A}} (\widehat{\boldsymbol{\beta}}(\text{AdaEnet})_{\mathcal{A}} - \boldsymbol{\beta}^*_{\mathcal{A}}).$$

Note that

$$\boldsymbol{\alpha}^T (\mathbf{I} + \lambda_2 \boldsymbol{\Sigma}^{-1}_{\mathcal{A}}) \boldsymbol{\Sigma}^{1/2}_{\mathcal{A}} \left(\widetilde{\boldsymbol{\beta}}^*_{\mathcal{A}} - \frac{\boldsymbol{\beta}^*_{\mathcal{A}}}{1 + \lambda_2/n}\right)$$

$$= \boldsymbol{\alpha}^T (\mathbf{I} + \lambda_2 \boldsymbol{\Sigma}^{-1}_{\mathcal{A}}) \boldsymbol{\Sigma}^{1/2}_{\mathcal{A}} \frac{\lambda_2 \boldsymbol{\beta}^*_{\mathcal{A}}}{n + \lambda_2} + \boldsymbol{\alpha}^T (\mathbf{I} + \lambda_2 \boldsymbol{\Sigma}^{-1}_{\mathcal{A}}) \boldsymbol{\Sigma}^{1/2}_{\mathcal{A}} (\widetilde{\boldsymbol{\beta}}^*_{\mathcal{A}} - \widetilde{\boldsymbol{\beta}}^*_{\mathcal{A}}(\lambda_2, 0))$$

$$+ \boldsymbol{\alpha}^T (\mathbf{I} + \lambda_2 \boldsymbol{\Sigma}^{-1}_{\mathcal{A}}) \boldsymbol{\Sigma}^{1/2}_{\mathcal{A}} (\widetilde{\boldsymbol{\beta}}^*_{\mathcal{A}}(\lambda_2, 0) - \boldsymbol{\beta}^*_{\mathcal{A}}).$$

In addition, we have

$$(\mathbf{I}+\lambda_2\mathbf{\Sigma}_{\mathcal{A}}^{-1})\mathbf{\Sigma}_{\mathcal{A}}^{1/2}(\widetilde{\boldsymbol{\beta}}_{\mathcal{A}}^*(\lambda_2,0)-\boldsymbol{\beta}_{\mathcal{A}}^*)=-\lambda_2\mathbf{\Sigma}_{\mathcal{A}}^{-1/2}\boldsymbol{\beta}_{\mathcal{A}}^*+\mathbf{\Sigma}_{\mathcal{A}}^{-1/2}\mathbf{X}_{\mathcal{A}}^T\boldsymbol{\varepsilon}.$$

Therefore, by Theorem 3.2, it follows that, with probability tending to 1, $z_n = T_1 + T_2 + T_3$, where

$$T_1 = \boldsymbol{\alpha}^T(\mathbf{I}+\lambda_2\mathbf{\Sigma}_{\mathcal{A}}^{-1})\mathbf{\Sigma}_{\mathcal{A}}^{1/2}\frac{\lambda_2\boldsymbol{\beta}_{\mathcal{A}}^*}{n+\lambda_2}-\boldsymbol{\alpha}^T\lambda_2\mathbf{\Sigma}_{\mathcal{A}}^{-1/2}\boldsymbol{\beta}_{\mathcal{A}}^*,$$

$$T_2 = \boldsymbol{\alpha}^T(\mathbf{I}+\lambda_2\mathbf{\Sigma}_{\mathcal{A}}^{-1})\mathbf{\Sigma}_{\mathcal{A}}^{1/2}(\widetilde{\boldsymbol{\beta}}_{\mathcal{A}}^*-\widetilde{\boldsymbol{\beta}}_{\mathcal{A}}^*(\lambda_2,0)),$$

$$T_3 = \boldsymbol{\alpha}^T\mathbf{\Sigma}_{\mathcal{A}}^{-1/2}\mathbf{X}_{\mathcal{A}}^T\boldsymbol{\varepsilon}.$$

We now show that $T_1 = o(1)$, $T_2 = o_P(1)$ and $T_3 \to N(0,\sigma^2)$ in distribution. Then, by Slutsky's theorem, we know $z_n \to_d N(0,\sigma^2)$. By (A1) and $\boldsymbol{\alpha}^T\boldsymbol{\alpha}=1$, we have

$$T_1^2 \leq 2\left\|(\mathbf{I}+\lambda_2\mathbf{\Sigma}_{\mathcal{A}}^{-1})\mathbf{\Sigma}_{\mathcal{A}}^{1/2}\frac{\lambda_2\boldsymbol{\beta}_{\mathcal{A}}^*}{n+\lambda_2}\right\|_2^2 + 2\|\lambda_2\mathbf{\Sigma}_{\mathcal{A}}^{-1/2}\boldsymbol{\beta}_{\mathcal{A}}^*\|_2^2$$

$$\leq 2\frac{\lambda_2^2}{(n+\lambda_2)^2}\|\mathbf{\Sigma}_{\mathcal{A}}^{1/2}\boldsymbol{\beta}_{\mathcal{A}}^*\|_2^2\left(1+\frac{\lambda_2}{bn}\right)^2 + 2\lambda^2\|\boldsymbol{\beta}_{\mathcal{A}}^*\|_2^2\frac{1}{bn}$$

$$\leq \frac{2\lambda_2^2 Bn}{(n+\lambda_2)^2}\left(1+\frac{\lambda_2}{bn}\right)^2\|\boldsymbol{\beta}_{\mathcal{A}}^*\|_2^2 + 2\lambda^2\|\boldsymbol{\beta}_{\mathcal{A}}^*\|_2^2\frac{1}{bn}.$$

Hence, it follows by (A6) that $T_1 = o(1)$. Similarly, we can bound $T_2$ as follows:

$$T_2^2 \leq \left(1+\frac{\lambda_2}{bn}\right)^2\|\mathbf{\Sigma}_{\mathcal{A}}^{1/2}(\widetilde{\boldsymbol{\beta}}_{\mathcal{A}}^*-\widetilde{\boldsymbol{\beta}}_{\mathcal{A}}^*(\lambda_2,0))\|_2^2$$

$$\leq \left(1+\frac{\lambda_2}{bn}\right)^2 Bn\|\widetilde{\boldsymbol{\beta}}_{\mathcal{A}}^*-\widetilde{\boldsymbol{\beta}}_{\mathcal{A}}^*(\lambda_2,0)\|_2^2$$

$$\leq \left(1+\frac{\lambda_2}{bn}\right)^2 Bn\left(\frac{\lambda_1^*\hat{\eta}^{-\gamma}}{bn+\lambda_2}\right)^2,$$

where we have used (6.10) in the last step. Then, (6.13) tells us that $T_2^2 = \frac{1}{n^2}O_P(1)$. Next, we consider $T_3$. Let $\mathbf{X}_{\mathcal{A}}[i,]$ denote the $i$th row of the matrix $\mathbf{X}_{\mathcal{A}}$. With such notation, we can write $T_3 = \sum_{i=1}^n r_i\varepsilon_i$, where $r_i = \boldsymbol{\alpha}^T(\mathbf{X}_{\mathcal{A}}^T\mathbf{X}_{\mathcal{A}})^{-1/2}(\mathbf{X}_{\mathcal{A}}[i,])^T$. Then, it is easy to see that

$$\sum_{i=1}^n r_i^2 = \sum_{i=1}^n \boldsymbol{\alpha}^T(\mathbf{X}_{\mathcal{A}}^T\mathbf{X}_{\mathcal{A}})^{-1/2}(\mathbf{X}_{\mathcal{A}}[i,])^T(\mathbf{X}_{\mathcal{A}}[i,])(\mathbf{X}_{\mathcal{A}}^T\mathbf{X}_{\mathcal{A}})^{-1/2}\boldsymbol{\alpha}$$

(6.14)
$$= \boldsymbol{\alpha}^T(\mathbf{X}_{\mathcal{A}}^T\mathbf{X}_{\mathcal{A}})^{-1/2}(\mathbf{X}_{\mathcal{A}}^T\mathbf{X}_{\mathcal{A}})(\mathbf{X}_{\mathcal{A}}^T\mathbf{X}_{\mathcal{A}})^{-1/2}\boldsymbol{\alpha}$$
$$= \boldsymbol{\alpha}^T\boldsymbol{\alpha} = 1.$$

Furthermore, we have for $k = 2 + \delta, \delta > 0$

$$\sum_{i=1}^n E[|\varepsilon_i|^{2+\delta}]|r_i^{2+\delta}| \leq E[|\varepsilon|^{2+\delta}]\left(\sum_{i=1}^n |r_i^2|\left(\max_i |r_i|^\delta\right)\right)$$

$$= E[|\varepsilon|^{2+\delta}]\left(\max_i |r_i^2|\right)^{\delta/2}.$$

Note that $r_i^2 \leq \|\mathbf{\Sigma}_{\mathcal{A}}^{-1/2}(\mathbf{X}_{\mathcal{A}}[i,])^T \leq (\sum_{j \in \mathcal{A}} x_{ij}^2)(\lambda_{\max}(\mathbf{\Sigma}_{\mathcal{A}}^{-1})) \leq \frac{\sum_{j=1}^p x_{ij}^2}{bn}$. Hence,

(6.15) $$\sum_{i=1}^n E[|\varepsilon_i|^{2+\delta}]|r_i^{2+\delta}| \leq E[|\varepsilon|^{2+\delta}]\left(\frac{\max_i(\sum_{j=1}^p x_{ij}^2)}{bn}\right)^{\delta/2} \to 0.$$

From (6.14) and (6.15), Lyapunov conditions for the central limit theorem are established. Thus, $T_3 \to_d N(0, \sigma^2)$. This completes the proof. $\square$

**Acknowledgments.** We sincerely thank an associate editor and referees for their helpful comments and suggestions.

SCHOOL OF STATISTICS  
UNIVERSITY OF MINNESOTA  
MINNEAPOLIS, MINNESOTA 55455  
USA  
E-MAIL: hzou@stat.umn.edu

DEPARTMENT OF STATISTICS  
NORTH CAROLINA STATE UNIVERSITY  
RALEIGH, NORTH CAROLINA 27695-8203  
USA  
E-MAIL: hzhang2@stat.ncsu.edu